\documentclass[11pt]{article}
\usepackage[utf8]{inputenc}
\usepackage[T1]{fontenc}
\usepackage[margin=1in]{geometry}
\usepackage{amsmath,amssymb,amsthm,mathtools}
\usepackage[dvipsnames]{xcolor}
\definecolor{citegreen}{RGB}{0,195,0}
\usepackage[colorlinks=true,
            linkcolor=blue,
            citecolor=citegreen,
            urlcolor=citegreen]{hyperref}
\usepackage[capitalise,nameinlink,noabbrev]{cleveref}
\usepackage[shortlabels]{enumitem}
\usepackage{tikz}
\usetikzlibrary{calc}
\usepackage{graphicx}

\theoremstyle{plain}
\newtheorem{theorem}{Theorem}
\newtheorem{lemma}[theorem]{Lemma}
\newtheorem{proposition}[theorem]{Proposition}
\newtheorem{corollary}[theorem]{Corollary}

\newtheorem{problem}[theorem]{Problem}

\newtheorem{claim}{Claim}[theorem]   

\newenvironment{claimproof}[1][Proof of the claim]
  {\begin{proof}[#1]}
  {\end{proof}}

\theoremstyle{definition}

\theoremstyle{remark}

\let\emptyset\varnothing

\newcommand{\ex}{\operatorname{ex}}
\newcommand{\ehom}{e_{\mathrm{hom}}}
\newcommand{\compl}[1]{#1^{\mathrm{c}}}

\title{A note on Tur\'an densities of daisies and hypergraph products}
\author{Marko Peji\'c\thanks{Karlsruhe Institute of Technology, Germany
  (\texttt{marko.pejic@student.kit.edu}).}}
\date{}

\begin{document}
\maketitle

\begin{abstract}
The Tur\'an density of an $r$-uniform hypergraph $\mathcal H$ is defined as
$\pi(\mathcal H)=\lim_{n\rightarrow\infty}\ex(n,\mathcal H)/\binom{n}{r}$,
where $\ex(n,\mathcal H)$ denotes the largest number of edges in an
$\mathcal H$-free $r$-uniform hypergraph on $n$ vertices.
For integers $r\geq s\geq2$ and $t\geq s$, an $(r,s,t)$-daisy
$\mathcal D_r(s,t)$ is the $r$-uniform hypergraph consisting of all
unions of a fixed $(r-s)$-set with the $s$-subsets of a disjoint
$t$-set.
We show that for every fixed $s\geq3$,
\[
    \lim_{r\rightarrow\infty}\pi(\mathcal D_r(s,t))=1-O_s(t^{-(s-1)})
\] 
as $t\rightarrow\infty$.
This bound is best possible up to the implicit constant by a
classical result of de~Caen, and it answers, in a stronger form,
a question of Ellis, Ivan and Leader, who asked for the case $s=3$. 
Moreover, we show a related result for matroid basis hypergraphs. 

We also consider Tur\'an densities of hypergraph products, of which daisies are a special case. 
The product of two hypergraphs on disjoint vertex sets consists
of all unions of an edge of the first with an edge of the second.
For the product $\mathcal H^d$ of $d$ disjoint copies of a fixed
$r$-uniform hypergraph $\mathcal H$ that is not $r$-partite, we show
that $1-\pi(\mathcal H^d) \rightarrow 0$ exponentially fast as
$d\rightarrow\infty$, and we determine the exponential rate in terms
of the minimum number of edges in a homomorphic image of $\mathcal H$.
This answers a question of Bollob\'as, Leader and Malvenuto up to a
polynomial factor in $d$, and we further address a related question of theirs
for products of two different hypergraphs.
\end{abstract}

\section{Introduction}

An \emph{$r$-uniform hypergraph} $\mathcal{H}$, or \emph{$r$-graph} for short,
consists of a finite vertex set $V(\mathcal{H})$ and an edge set
$E(\mathcal{H})\subseteq\binom{V(\mathcal{H})}{r}$, where
$\binom{X}{k}$ denotes the family of all subsets of $X$ of size $k$. 
The complete $r$-graph on $n$ vertices is denoted by $K_{n}^{(r)}$,
and we may also refer to $2$-graphs simply as graphs. 
For an $r$-graph $\mathcal{G}$ and a set $S\subseteq V(\mathcal{G})$
with $|S|<r$, the \emph{link} of $S$ is the $(r-|S|)$-graph on
$V(\mathcal{G})\setminus S$ with edge set
$\{B\subseteq V(\mathcal{G})\setminus S: S\cup B\in E(\mathcal{G})\}$.
We call links of single vertices \emph{vertex links}. 
For $r$-graphs $\mathcal{G}$ and $\mathcal{H}$, a \emph{copy} of
$\mathcal{H}$ in $\mathcal{G}$ is a subhypergraph of $\mathcal{G}$
isomorphic to $\mathcal{H}$, and $\mathcal{G}$ is
\emph{$\mathcal{H}$-free} if it contains no copy of $\mathcal{H}$.
The \emph{extremal number} $\ex(n,\mathcal{H})$ is the largest number
of edges in an $\mathcal{H}$-free $r$-graph on $n$ vertices.
The \emph{Tur\'an density} of $\mathcal{H}$ is defined as
\[
    \pi(\mathcal{H})=\lim_{n\rightarrow\infty}\frac{\ex(n,\mathcal{H})}{\binom{n}{r}}.
\]
Here and throughout, the uniformity is always that of the hypergraph
in question, and we do not indicate it explicitly in the notation.
By a standard double counting argument, the sequence
$\ex(n,\mathcal{H})/\binom{n}{r}$ is nonincreasing in $n$, so this density
notion is well defined.
We assume that all hypergraphs whose Tur\'an density we consider have
at least one edge and no isolated vertices.
In asymptotic notation, the implicit constants may depend on the
parameters given in the subscript. 
Throughout, $\log$ denotes the natural logarithm.

\paragraph{Daisies.}
For integers $r\geq s\geq2$ and $t\geq s$, an \emph{$(r,s,t)$-daisy}
is an $r$-graph with vertex set $S\cup T$ and edge set
$\{S\cup X: X\in\binom{T}{s}\}$, where $S\cap T=\emptyset$,
$|S|=r-s$, and $|T|=t$.
We denote this hypergraph by $\mathcal D_r(s,t)$\footnote{We follow
the notation of Ellis, Ivan and Leader~\cite{EIL}. In the original
work of Bollob\'as, Leader and Malvenuto~\cite{BLM}, which introduced
daisies, the roles of $s$ and $t$ are interchanged, that is, their
$(s,t)$-daisy is our $\mathcal{D}_r(t,s)$.} and call $S$ its
\emph{stem}.
The daisy $\mathcal D_r(s,t)$ has $r-s+t$ vertices and $\binom{t}{s}$
edges.
A particular question for daisies, raised by Bollob\'as, Leader and
Malvenuto~\cite{BLM}, is how their Tur\'an densities behave as the
size of the stem grows.
For fixed $s$ and $t$, averaging over vertex links shows that
$\pi(\mathcal D_r(s,t))$ is nonincreasing in $r$, and thus
$\lim_{r\rightarrow\infty}\pi(\mathcal D_r(s,t))$ exists.

It was conjectured in~\cite{BLM}, and independently by Bukh as
mentioned there, that this limit is zero for every fixed pair $s$ and
$t$ with $t\geq s$.
Using a construction based on bases of an $r$-dimensional vector space over a finite field, Ellis, Ivan and Leader~\cite{EIL}
disproved this conjecture by showing that
$\lim_{r\to\infty}\pi(\mathcal D_r(2,4))
\geq\prod_{k=1}^{\infty}(1-2^{-k})\approx0.29$.
In the same work, they further showed that
$\lim_{r\rightarrow\infty}\pi(\mathcal D_r(2,t))
=1-1/t+o(1/t)$ as $t\rightarrow\infty$.
The upper bound $\pi(\mathcal D_r(2,t))\leq1-1/(t-1)$ used for this result
follows by an application of Tur\'an's theorem~\cite{Tur} to the links of
$(r-2)$-sets and was recently improved by Hou, Hou, Liu, Zeng and
Zhang~\cite{HHLZZ} to
$\pi(\mathcal D_r(2,t))
\leq1-1/(t-1)-1/(12(t-1)^2)$ for all $r\geq3$ and $t\geq3$.
For larger even values of $s$, Ellis, Ivan and Leader also obtained
lower bounds by considering vector spaces of dimension larger than $r$.

The underlying idea of their construction has since found further
applications, for example to Tur\'an problems in the
hypercube, see~\cite{AP,EIL2,GM,Pej}, while the construction itself has
also been used in the context of erasure list-decodable codes by
Alon~\cite{Alon} and of matroid basis hypergraphs by van der Pol,
Walsh and Wigal~\cite{PWW}.
More recently, Pikhurko, Yang and Zhu~\cite{PYZ} considered
hypergraph augmentations, which are obtained by adding the same set of
new vertices to every edge of a fixed hypergraph.
Their work implies, in particular, that
$\pi(\mathcal D_r(s,s+1))=\Theta_s(1/r)$ as $r\rightarrow\infty$
for every fixed $s\geq2$.

However, for fixed $s\geq3$, the aforementioned results do not
determine the correct order of
$1-\lim_{r\rightarrow\infty}\pi(\mathcal D_r(s,t))$
as $t\rightarrow\infty$.
In fact, Ellis, Ivan and Leader~\cite{EIL} explicitly asked the
following question for $s=3$.

\begin{problem}\label{problem: daisy}
Is it true that $\lim_{r\rightarrow\infty}\pi(\mathcal D_r(3,t))=1-O(t^{-2})$ as $t\rightarrow\infty$?
\end{problem}

Our first goal in this note is to answer this question, and in fact,
to do so for every fixed $s\geq3$. 
More specifically, we prove the following result, whose case $s=3$
answers Problem~\ref{problem: daisy} in the affirmative.

\begin{theorem}\label{main}
For every fixed integer $s\geq3$,
\[
    \lim_{r \rightarrow \infty}\pi(\mathcal{D}_r(s,t)) = 1 - O_s\!\left(t^{-(s-1)}\right)
\]
as $t \rightarrow \infty$.
\end{theorem}

Note that, by monotonicity in $r$, the same lower bound holds uniformly
for all $r\geq s$.
Also, this bound is best possible up to the implicit constant,
since the bound of de~Caen~\cite{dC} for $K^{(s)}_t$ and averaging
over the links of $(r-s)$-sets give
$\pi(\mathcal{D}_r(s,t))\leq1-\binom{t-1}{s-1}^{-1}$
for all $r\geq s$, as discussed for $s=3$ in~\cite{EIL}. 
We shall also see how the ideas used to prove Theorem~\ref{main}
can be applied to matroid basis hypergraphs.

\paragraph{Hypergraph products.}
A more general problem related to daisies, concerning \emph{hypergraph products}, was raised by Bollob\'as, Leader and Malvenuto~\cite{BLM}.
Given an $r$-graph $\mathcal{F}$ and an $s$-graph $\mathcal{G}$ on
disjoint vertex sets, their \emph{product} $\mathcal{F}*\mathcal{G}$ is
the $(r+s)$-graph on $V(\mathcal{F})\cup V(\mathcal{G})$ with edge set
\[
E(\mathcal{F}*\mathcal{G})=\{A\cup B: A\in E(\mathcal{F}),\ B\in E(\mathcal{G})\}.
\]
Note that a daisy is itself a product, since for $r>s$ the daisy
$\mathcal{D}_r(s,t)$ is the product of $K_{t}^{(s)}$ with the
$(r-s)$-graph consisting of a single edge.
Bollob\'as, Leader and Malvenuto~\cite[Problem~5]{BLM} asked the
following.

\begin{problem}\label{problem: two_product}
Let $\mathcal{F}$ be an $r$-graph and $\mathcal{G}$ be an $s$-graph.
How does $\pi(\mathcal{F}*\mathcal{G})$ compare to $\pi(\mathcal{F})$
and $\pi(\mathcal{G})$?
\end{problem}

A particularly interesting case, as pointed out in~\cite{BLM}, is that
of products of copies of a single hypergraph.
Writing $\mathcal{F}^d$ for the product of $d$ disjoint copies of
$\mathcal{F}$, which we shall call the $d$-th \emph{power} of $\mathcal{F}$,
Bollob\'as, Leader and Malvenuto~\cite[Problem~6]{BLM} also asked the following.

\begin{problem}\label{problem: fold}
Let $\mathcal{F}$ be a fixed $r$-graph. As $d$ varies, how does
$\pi(\mathcal{F}^d)$ behave? In particular, what happens when
$\mathcal{F}=K_{t}^{(r)}$?
\end{problem}

Hypergraph products have appeared in the context of Betti
numbers~\cite{Emtander} and Lagrangians~\cite{Grosu},
and results on their Tur\'an densities are also known in the special
case where one factor consists of a single
edge, see~\cite{FRV,EIL,CGKZ,PYZ}.
However, to the best of our knowledge,
Problems~\ref{problem: two_product} and~\ref{problem: fold}
have not previously been addressed as originally posed.
Our second goal in this note is to make progress on them,
with the following parameter as the key ingredient.

A \emph{homomorphism} from an $r$-graph $\mathcal{F}$ to an $r$-graph
$\mathcal{H}$ is a map $\varphi:V(\mathcal{F})\rightarrow V(\mathcal{H})$
such that $\varphi(A)\in E(\mathcal{H})$ for every $A\in E(\mathcal{F})$. 
In particular, $\varphi$ is injective on every edge. 
We write $\mathcal{F}\rightarrow\mathcal{H}$ if a homomorphism from
$\mathcal{F}$ to $\mathcal{H}$ exists, and define
\[
    \ehom(\mathcal{F}) =\min\{|E(\mathcal{H})|:\mathcal{F}\rightarrow\mathcal{H}\}
\]
as the \emph{homomorphic edge number} of $\mathcal{F}$.
In other words, $\ehom(\mathcal{F})$ denotes the least number of edges in a homomorphic image of $\mathcal{F}$. 
For example, $\ehom(\mathcal{F})=1$ if and only if $\mathcal{F}$ is
$r$-partite.
Another example is given by $\ehom(C_{2\ell+1})=3$ for every $\ell\geq1$, since an odd cycle maps onto a triangle, and every homomorphic image of an odd
cycle is non-bipartite and therefore has at least three edges. 
We assume that all hypergraphs whose homomorphic edge number we consider have at least one edge and no isolated vertices.
With this parameter at hand, we prove the following result on Problem~\ref{problem: fold}.

\begin{theorem}\label{theorem: products}
Let $\mathcal{F}$ be an $r$-graph.
If $\mathcal{F}$ is $r$-partite, then $\pi(\mathcal{F}^d)=0$ for all
$d\geq1$.
Otherwise, there is a constant $C=C(\mathcal{F})>0$ such that, for every
integer $d\geq1$,
\[
    \ehom(\mathcal{F})^{-d}
    \leq 1-\pi(\mathcal{F}^d)
    \leq Cd^2\,\ehom(\mathcal{F})^{-d}.
\]
\end{theorem}

The lower bound follows from a standard averaging argument
applied to a homomorphic image of $\mathcal{F}^d$.
Our main contribution here is to show that it is sharp up to a
polynomial factor in $d$ and a constant $C$ depending only on
$\mathcal{F}$.

Observe that $\ehom(K^{(r)}_t)=\binom{t}{r}$ for $t\geq r\geq2$,
since every pair of vertices of $K^{(r)}_t$ is contained in an
edge and thus every homomorphism from $K^{(r)}_t$ is injective.
Also, $\ehom(G)=\binom{\chi(G)}{2}$ for every graph $G$, since a proper vertex-colouring with $\chi(G)$ colours is a
homomorphism $G\rightarrow K_{\chi(G)}$, and conversely every
homomorphic image $H$ of $G$ satisfies $\chi(H)\geq\chi(G)$.
Since any two colour classes in a proper vertex-colouring of $H$
with $\chi(H)$ colours are joined by an edge, it follows that
$|E(H)|\geq\binom{\chi(H)}{2}\geq\binom{\chi(G)}{2}$.
Using these observations, the following two results are immediate
consequences of Theorem~\ref{theorem: products}.

\begin{corollary}\label{corollary: cliques}
Let $t>r\geq2$ be integers. There is a constant $C=C(r,t)>0$ such that,
for every integer $d\geq1$,
\[
    \binom{t}{r}^{-d}
    \leq 1-\pi\bigl((K^{(r)}_t)^d\bigr)
    \leq Cd^2\binom{t}{r}^{-d}.
\]
\end{corollary}

\begin{corollary}\label{corollary: graphs}
Let $G$ be a graph. If $G$ is bipartite, then $\pi(G^d)=0$ for all
$d\geq1$. Otherwise, there is a constant $C=C(G)>0$ such that, for
every integer $d\geq1$,
\[
    \binom{\chi(G)}{2}^{-d}
    \leq 1-\pi(G^d)
    \leq Cd^2\binom{\chi(G)}{2}^{-d}.
\]
\end{corollary}

For Problem~\ref{problem: two_product}, we only make partial progress
in the form of the following result, which, however, also points to
the homomorphic edge number as a relevant parameter.

\begin{proposition}\label{theorem: two factors}
Let $\mathcal{F}$ be an $r$-graph and $\mathcal{G}$ be an $s$-graph.
Then
\[
    1-\pi(\mathcal{F}*\mathcal{G})
    \geq\max\left\{
        \frac{1-\pi(\mathcal{G})}{\ehom(\mathcal{F})},\
        \frac{1-\pi(\mathcal{F})}{\ehom(\mathcal{G})}
    \right\}.
\]
\end{proposition}

In the concluding remarks, we shall also give an example, based on
Proposition~\ref{theorem: two factors}, showing that
$\pi(\mathcal{F}*\mathcal{G})$ is not determined by $\pi(\mathcal{F})$,
$\pi(\mathcal{G})$ and the uniformities of $\mathcal{F}$ and
$\mathcal{G}$ alone.

\paragraph{Organisation.}
The remainder of this note is structured as follows.
In Section~\ref{section: theorem_daisy}, we prove Theorem~\ref{main},
and in Section~\ref{subsection: matroids}, we show a related result
for matroid basis hypergraphs.
In Section~\ref{section: products}, we prove
Theorem~\ref{theorem: products} and Proposition~\ref{theorem: two factors}, after
giving two auxiliary results on homomorphisms.
We close this note with some concluding remarks and further directions
in Section~\ref{section: remarks}.

\section{Daisies}\label{section: daisy}

The main purpose of this section is to prove Theorem~\ref{main}.
Our construction is essentially the same as that used by Ellis,
Ivan and Leader~\cite{EIL} to prove their Theorems~5 and~10,
with only a slight modification.
The difference is that instead of retaining linearly independent sets
of $r$ vectors in a space of dimension at least $r$, we retain sets of
$r$ vectors \emph{spanning} $\mathbb F_q^{r-s+2}$.
The reason for this is that spanning a space of dimension $r-s+2$ is a
much weaker requirement on $r$ vectors than linear independence, and
thus more $r$-sets are retained.
More specifically, $r$ vectors chosen independently and uniformly at
random from $\mathbb F_q^{r}$ fail to be linearly independent with
probability about $1/q$, whereas $r$ vectors chosen in the same way
from $\mathbb F_q^{r-s+2}$ fail to span it with probability only about
$q^{-(s-1)}$, which also explains where the exponent $s-1$ comes from.
Also, the smaller the dimension of the space, the more likely it is
that $r$ random vectors span it, so we want the dimension to be as small as possible. 
At the same time, the dimension has to be large enough to avoid copies
of $\mathcal{D}_r(s,t)$, and for our construction this means that the
quotient by the span of the vectors of a stem has dimension at least
two.
As a stem has $r-s$ vertices, this gives the dimension $r-s+2$.

We make this precise in Section~\ref{section: theorem_daisy}, and
in Section~\ref{subsection: matroids} we apply the construction to
matroid basis hypergraphs.

\subsection{Proof of \texorpdfstring{Theorem~\ref{main}}{Theorem 2}}\label{section: theorem_daisy}

We begin with the following proposition, which is the main ingredient
in the proof of Theorem~\ref{main}.

\begin{proposition}\label{proposition: construction}
Let $r\geq s\geq2$ be integers. Let $q$ be a prime power and $t$ an
integer such that $t>(s-1)(q+1)$.
For every integer $n\geq r$,
\[
    \ex(n,\mathcal D_r(s,t))
    \geq
    \left(1-\frac{1}{q^{s-2}(q-1)}\right)\binom{n}{r}.
\]
\end{proposition}

\begin{proof}[Proof of Proposition~\ref{proposition: construction}]
Let $d=r-s+2$.
For every $i\in[n]=\{1,\ldots,n\}$, pick a vector
$v_i\in\mathbb F_q^d$, independently and uniformly at random.
Retain an $r$-set $A\subseteq[n]$ if
$\operatorname{span}\{v_i:i\in A\}=\mathbb F_q^d$, and let
$\mathcal{G}$ be the $r$-graph on $[n]$ whose edges are the retained
sets.
Note that the vectors are not required to be distinct or nonzero, and
that $\mathcal{G}$ depends on the choice of the vectors.

We also note that the random choice of the vectors is only needed to show that
$\mathcal{G}$ has many edges for some choice, whereas $\mathcal{G}$
contains no copy of $\mathcal{D}_r(s,t)$ for every choice.
We prove these two properties in the following two claims.

\begin{claim}\label{claim: one}
$\mathcal{G}$ contains no copy of $\mathcal D_r(s,t)$ for every choice of the vectors.
\end{claim}
\begin{claimproof}[Proof of Claim~\ref{claim: one}]
Suppose otherwise, and let
$\{S\cup X:X\in\binom{T}{s}\}$ be such a copy, where
$S\cap T=\emptyset$, $|S|=r-s$, and $|T|=t$. 
Since $|S|=d-2$, there is a subspace $W\subseteq\mathbb F_q^d$ of
dimension $d-2$ containing all vectors $v_i$ with $i\in S$.
Let $U=\mathbb F_q^d/W$ and let
$\varphi:\mathbb F_q^d\rightarrow U$ be the natural projection onto $U$, given by $\varphi(x)=x+W$ for
$x\in\mathbb F_q^d$.
The space $U$ has dimension two and thus contains
$q^2-1$ nonzero vectors.
Each nonzero vector of $U$ is contained in exactly one subspace of
dimension one, namely its span, and each such subspace has $q-1$
nonzero vectors. 
Thus $U$ has $(q^2-1)/(q-1)=q+1$ subspaces of
dimension one. 
We define a colouring of $T$ with these $q+1$ subspaces as follows.
For each $i\in T$ with $\varphi(v_i)\neq0$, give $i$ the
colour of the unique subspace containing $\varphi(v_i)$.
If $\varphi(v_i)=0$, give $i$ any fixed colour.
Since $t > (s-1)(q+1)$, the pigeonhole principle gives a set
$X \in \binom{T}{s}$ whose elements all have the same colour.
Thus the vectors $\varphi(v_i)$ with $i\in X$ are contained in
$\operatorname{span}\{u\}$ for some nonzero $u\in U$.
Since $v_i\in W$ for every $i\in S$, we have
$v_i\in\varphi^{-1}(\operatorname{span}\{u\})$ for all $i\in S\cup X$. 
At the same time, since $\ker\varphi=W$, the subspace $\varphi^{-1}(\operatorname{span}\{u\})$
has dimension $\dim W+1=d-1$.
Thus $\dim\operatorname{span}\{v_i:i\in S\cup X\}\leq d-1$, so the
vectors $v_i$ with $i\in S\cup X$ do not span $\mathbb F_q^d$, a
contradiction, since $S\cup X\in E(\mathcal{G})$.
This proves Claim~\ref{claim: one}.
\end{claimproof}

\begin{claim}\label{claim: two}
There exists a choice of the vectors such that $|E(\mathcal{G})| \geq \left(1-1/({q^{s-2}(q-1)})\right)\binom{n}{r}$.
\end{claim}

\begin{claimproof}[Proof of Claim~\ref{claim: two}]
Consider an $r$-set $A\subseteq[n]$.
If $A$ is not retained, then all vectors $v_i$ with $i\in A$ are
contained in a subspace of dimension $d-1$.
We claim that there are $(q^d-1)/(q-1)$ such subspaces.
Indeed, a linear map $f:\mathbb F_q^d\rightarrow\mathbb F_q$ is given
by $f(x_1,\ldots,x_d)=\lambda_1x_1+\cdots+\lambda_dx_d$,
where $\lambda_1,\ldots,\lambda_d\in\mathbb F_q$, and its kernel is the subspace
$\ker f =\{x\in\mathbb F_q^d: f(x)=0\}$.
If $f$ is nonzero, then $\ker f$ has dimension $d-1$, and conversely
every subspace of dimension $d-1$ is the kernel of a nonzero linear
map.
There are $q^d-1$ such maps, and two have the same kernel exactly
when one is a nonzero scalar multiple of the other.
Thus each subspace is counted $q-1$ times, giving $(q^d-1)/(q-1)$
subspaces, and each of these subspaces contains $q^{d-1}$ vectors.
Since the vectors are chosen independently and uniformly,
the probability that all $v_i$ with $i\in A$ are contained in
any fixed one of these subspaces is $q^{-r}$.
By the union bound,
\[
    \Pr(A\notin E(\mathcal{G}))
    \leq\frac{q^d-1}{q-1}q^{-r}
    \leq\frac{1}{q^{s-2}(q-1)},
\]
where the last inequality uses $d=r-s+2$.
Thus, by linearity of expectation,
\[
    \mathbb E(|E(\mathcal{G})|)
    \geq
    \left(1-\frac{1}{q^{s-2}(q-1)}\right)\binom{n}{r},
\]
so some choice of the vectors gives at least this many edges. This proves Claim~\ref{claim: two}.
\end{claimproof}
Together, Claims~\ref{claim: one} and~\ref{claim: two} prove Proposition~\ref{proposition: construction}.
\end{proof}

\paragraph{Remark.} For $s=2$, the vector space $\mathbb F_q^{r-s+2}$ has dimension
$r$, so the retained sets are exactly the linearly independent ones
and Proposition~\ref{proposition: construction} reduces to the construction of Ellis, Ivan and Leader, with the slightly weaker bound $1-1/(q-1)$.

\begin{proof}[Proof of Theorem~\ref{main}]
Fix $s\geq3$ and let $t\geq3s-2$.
Pick $q$ to be a power of two with $q\leq(t-s)/(s-1)<2q$.
Then $q\geq2$, $(s-1)(q+1)\leq t-1<t$, and $q=\Theta_s(t)$.
Proposition~\ref{proposition: construction} therefore applies
for every $r\geq s$ and $n\geq r$.
Letting $n\rightarrow\infty$ and then $r\rightarrow\infty$,
with $s$, $t$, and $q$ fixed, gives
\[
    1-\lim_{r\rightarrow\infty}\pi(\mathcal D_r(s,t))
    \leq\frac{1}{q^{s-2}(q-1)}
    \leq\frac{2}{q^{s-1}}
    =O_s\!\left(t^{-(s-1)}\right),
\]
which implies the claimed bound.
\end{proof}

\subsection{An application to matroid basis hypergraphs}\label{subsection: matroids}

Our construction from Section~\ref{section: theorem_daisy} also
gives a bound for matroid basis hypergraphs. 
For a rank-$r$ matroid $M$, its \emph{basis hypergraph} is the
$r$-graph on the ground set of $M$ whose edges are the bases of $M$.
For a more detailed introduction to matroids and the standard terminology, we refer to Oxley~\cite{Oxley}.
Van der Pol, Walsh and Wigal~\cite{PWW} considered Tur\'an problems
for matroid basis hypergraphs and showed that the basis hypergraph
of a rank-$r$ matroid contains $\mathcal D_r(s,t)$ if and only if
the matroid has a $U_{s,t}$-minor, see~\cite[Corollary~2.3]{PWW}. 
Here, $U_{s,t}$ denotes the uniform matroid of rank $s$ on $t$ elements.

For integers $n\geq r\geq s\geq1$ and $t\geq s$, let
$\ex_{\mathsf{M}}(n,r,U_{s,t})$ be the largest number of bases in a
rank-$r$ matroid on $n$ elements containing no $U_{s,t}$-minor, where $\ex_{\mathsf{M}}(n,r,U_{s,t})=0$ if no such matroid exists.
Following the notation of van der Pol, Walsh and Wigal~\cite{PWW}, define
\[
    \pi_{\mathsf{M}}(r,U_{s,t})
    =\lim_{n\rightarrow\infty}
    \frac{\ex_{\mathsf{M}}(n,r,U_{s,t})}{\binom{n}{r}}.
\]
Both this limit and its limit as $r\rightarrow\infty$ exist, see
\cite[Corollary~3.3]{PWW}.
For positive integers $s$ and $t$ with $t\geq s$,
van der Pol, Walsh and Wigal asked for the value of
$\lim_{r\rightarrow\infty}\pi_{\mathsf{M}}(r,U_{s,t})$,
see~\cite[Problem~1.1(iii)]{PWW}.

Motivated by their connection between daisies and matroid basis
hypergraphs, we give an asymptotic answer to this question by
showing that our construction from Section~\ref{section: theorem_daisy}
also yields matroids with many bases and no $U_{s,t}$-minor.
More specifically, we determine how fast
$\lim_{r\rightarrow\infty}\pi_{\mathsf{M}}(r,U_{s,t})$
approaches one as $t\rightarrow\infty$, up to a factor
depending only on $s$, for every fixed $s\geq3$.

\begin{proposition}\label{corollary:matroid_basis}
For every fixed integer $s\geq3$,
\[
    \lim_{r\rightarrow\infty}\pi_{\mathsf{M}}(r,U_{s,t})
    =1-\Theta_s\!\left(t^{-(s-1)}\right)
\]
as $t\rightarrow\infty$.
\end{proposition}

\begin{proof}[Proof of Proposition~\ref{corollary:matroid_basis}]
Fix $s\geq3$ and $t\geq3s-2$, pick $q$ to be a power of two with
$q\leq(t-s)/(s-1)<2q$ as in the proof of Theorem~\ref{main}, and let
$n\geq r\geq s$.
By Claim~\ref{claim: two}, there are vectors
$v_1,\ldots,v_n\in\mathbb F_q^{r-s+2}$ for which the $r$-graph
$\mathcal{G}$ on $[n]$, whose edges are the $r$-sets $A\subseteq[n]$
with $\operatorname{span}\{v_i:i\in A\}=\mathbb F_q^{r-s+2}$, has at
least $\bigl(1-1/(q^{s-2}(q-1))\bigr)\binom{n}{r}$ edges.
We first show that its edges are the bases of a matroid.

To this end, let $A,B\in E(\mathcal{G})$ and $i\in A\setminus B$, and let
$W=\operatorname{span}\{v_j:j\in A\setminus\{i\}\}$.
If $W=\mathbb F_q^{r-s+2}$, then
$(A\setminus\{i\})\cup\{j\}\in E(\mathcal{G})$ for every
$j\in B\setminus A$.
Otherwise, $\dim W=r-s+1$, since the vectors indexed by $A$ span
$\mathbb F_q^{r-s+2}$ and adding one vector can increase the
dimension of a span by at most one.
Since the vectors $v_k$ with $k\in B$ span $\mathbb F_q^{r-s+2}$, there
is some $j\in B$ with $v_j\notin W$.
As $v_k\in W$ for all $k\in A\setminus\{i\}$ and $i\notin B$, we have
$j\in B\setminus A$.
Again, $(A\setminus\{i\})\cup\{j\}\in E(\mathcal{G})$.
Thus $E(\mathcal{G})$ is the set of bases of a rank-$r$ matroid $M$.

By~\cite[Corollary~2.3]{PWW}, a rank-$r$ matroid contains a
$U_{s,t}$-minor if and only if its basis hypergraph contains
$\mathcal D_r(s,t)$.
Since $\mathcal{G}$ is $\mathcal D_r(s,t)$-free by
Claim~\ref{claim: one}, $M$ has no $U_{s,t}$-minor, and thus
\[
    \ex_{\mathsf{M}}(n,r,U_{s,t})
    \geq
    \left(1-\frac{1}{q^{s-2}(q-1)}\right)\binom{n}{r}.
\]
Letting $n\rightarrow\infty$ and then $r\rightarrow\infty$, and using
$q=\Theta_s(t)$, gives
$1-\lim_{r\rightarrow\infty}\pi_{\mathsf{M}}(r,U_{s,t})=O_s\!\left(t^{-(s-1)}\right)$.

By the other direction of~\cite[Corollary~2.3]{PWW}, the basis
hypergraph of a rank-$r$ matroid without a $U_{s,t}$-minor is
$\mathcal D_r(s,t)$-free, so
$\pi_{\mathsf{M}}(r,U_{s,t})\leq\pi(\mathcal D_r(s,t))$.
The matching lower bound on
$1-\lim_{r\rightarrow\infty}\pi_{\mathsf{M}}(r,U_{s,t})$ now follows from
$\pi(\mathcal{D}_r(s,t))\leq1-\binom{t-1}{s-1}^{-1}$, which in turn follows from the result of de~Caen~\cite{dC}, as noted in the introduction.
\end{proof}

\section{Hypergraph products}\label{section: products}

We now prove Theorem~\ref{theorem: products}
and Proposition~\ref{theorem: two factors} in Sections~\ref{section: powers}
and~\ref{section: two factors}, respectively.
Before doing so, we give in Section~\ref{section: homomorphisms} two
lemmas on homomorphisms needed for these proofs.

\subsection{Auxiliary results}\label{section: homomorphisms}

The first of the two lemmas, Lemma~\ref{lemma: standard}, is the
standard observation that the Tur\'an density is monotone under
homomorphisms, for which we give a short proof for completeness, and
the second, Lemma~\ref{lemma: edge_number_product}, shows that the
homomorphic edge number is multiplicative with respect to products.

Here and in the remainder of this note, we shall also need the notion of a
\emph{blow-up} of an $r$-graph $\mathcal{H}$, which is obtained by
replacing every vertex of $\mathcal{H}$ by a set of vertices, called a
\emph{part}, and every edge of $\mathcal{H}$ by all $r$-sets meeting
each of the corresponding parts in exactly one vertex.
If every part has $t\geq1$ vertices, we denote the blow-up by
$\mathcal{H}(t)$.

\begin{lemma}\label{lemma: standard}
Let $\mathcal{F}$ and $\mathcal{H}$ be $r$-graphs.
If $\mathcal{F}\rightarrow\mathcal{H}$, then $\pi(\mathcal{F})\leq\pi(\mathcal{H})$.
Moreover, $\pi(\mathcal{F})=0$ if and only if $\mathcal{F}$ is $r$-partite.
\end{lemma}

\begin{proof}[Proof of Lemma~\ref{lemma: standard}] 
A standard supersaturation result shows that
$\pi(\mathcal{H}(t))=\pi(\mathcal{H})$ (see, e.g.,~\cite[Theorem~2.2]{Keevash}).
If $\mathcal{F}\rightarrow\mathcal{H}$, then $\mathcal{H}(t)$ contains
a copy of $\mathcal{F}$ for $t=|V(\mathcal{F})|$, obtained by mapping
vertices of $\mathcal{F}$ with the same image to distinct vertices of
the corresponding part. 
Thus $\pi(\mathcal{F})\leq\pi(\mathcal{H}(t))=\pi(\mathcal{H})$. 
Moreover, if $\mathcal{F}$ is $r$-partite, then
$\mathcal{F}\rightarrow K^{(r)}_r$, and thus
$\pi(\mathcal{F})\leq\pi(K^{(r)}_r)=0$.
If $\mathcal{F}$ is not $r$-partite, then complete $r$-partite
$r$-graphs show that $\pi(\mathcal{F})>0$.
\end{proof}

\begin{lemma}\label{lemma: edge_number_product}
Let $\mathcal{F}$ be an $r$-graph and $\mathcal{G}$ be an $s$-graph.
Then $\ehom(\mathcal{F}*\mathcal{G})=\ehom(\mathcal{F})\,\ehom(\mathcal{G})$.
In particular, $\ehom(\mathcal{F}^d)=\ehom(\mathcal{F})^d$ for every
integer $d\geq1$.
\end{lemma}

\begin{proof}[Proof of Lemma~\ref{lemma: edge_number_product}] 
We first show that
$\ehom(\mathcal{F}*\mathcal{G})\geq\ehom(\mathcal{F})\,\ehom(\mathcal{G})$. 
To this end, let $\varphi$ be a homomorphism from $\mathcal{F}*\mathcal{G}$ to an
$(r+s)$-graph $\mathcal{H}$ with
$|E(\mathcal{H})|=\ehom(\mathcal{F}*\mathcal{G})$.
By replacing $\mathcal{H}$ with the subhypergraph consisting of the
vertices and edges in the image of $\varphi$, we may assume that every
edge of $\mathcal{H}$ is the image of an edge of
$\mathcal{F}*\mathcal{G}$.
Moreover, since $\mathcal{F}$ and $\mathcal{G}$ have no isolated
vertices, every vertex of $\mathcal{F}$ is contained in a common edge
of $\mathcal{F}*\mathcal{G}$ with every vertex of $\mathcal{G}$, and
$\varphi$ is injective on every edge, so $\varphi(V(\mathcal{F}))$ and
$\varphi(V(\mathcal{G}))$ are disjoint.
Let $\mathcal{H}_1$ be the $r$-graph on $\varphi(V(\mathcal{F}))$ with
edge set $\{\varphi(A):A\in E(\mathcal{F})\}$, and let $\mathcal{H}_2$
be the $s$-graph on $\varphi(V(\mathcal{G}))$ with edge set
$\{\varphi(B):B\in E(\mathcal{G})\}$.
Then $\mathcal{H}_1$ and $\mathcal{H}_2$ are homomorphic images of
$\mathcal{F}$ and $\mathcal{G}$, respectively, and thus
$|E(\mathcal{H}_1)|\geq\ehom(\mathcal{F})$ and
$|E(\mathcal{H}_2)|\geq\ehom(\mathcal{G})$.
Every edge of $\mathcal{H}$ is of the form $\varphi(A)\cup\varphi(B)$
with $A\in E(\mathcal{F})$ and $B\in E(\mathcal{G})$, and conversely
every such set is an edge of $\mathcal{H}$, as $\varphi$ is a
homomorphism.
Since $V(\mathcal{H}_1)=\varphi(V(\mathcal{F}))$ and
$V(\mathcal{H}_2)=\varphi(V(\mathcal{G}))$ are disjoint, the set
$\varphi(A)\cup\varphi(B)$ determines both $\varphi(A)$ and
$\varphi(B)$.
Thus
\[
    \ehom(\mathcal{F}*\mathcal{G})
    =|E(\mathcal{H})|
    =|E(\mathcal{H}_1)| \cdot |E(\mathcal{H}_2)|
    \geq\ehom(\mathcal{F}) \cdot \ehom(\mathcal{G}).
\]

For the other direction, let $\mathcal{H}_1'$ and $\mathcal{H}_2'$ be
homomorphic images of $\mathcal{F}$ and $\mathcal{G}$ on disjoint
vertex sets with $\ehom(\mathcal{F})$ and $\ehom(\mathcal{G})$ edges,
respectively.
Then the two homomorphisms together give a homomorphism from
$\mathcal{F}*\mathcal{G}$ to $\mathcal{H}_1'*\mathcal{H}_2'$, which
has $\ehom(\mathcal{F})\,\ehom(\mathcal{G})$ edges.
Finally, the statement for $\mathcal{F}^d$ follows by induction on $d$.
\end{proof}

\subsection{Proof of \texorpdfstring{Theorem~\ref{theorem: products}}{Theorem 5}}\label{section: powers}

Since the idea behind the proof is easily lost among the
technicalities and the parameters involved, we first sketch the
argument before giving the details.

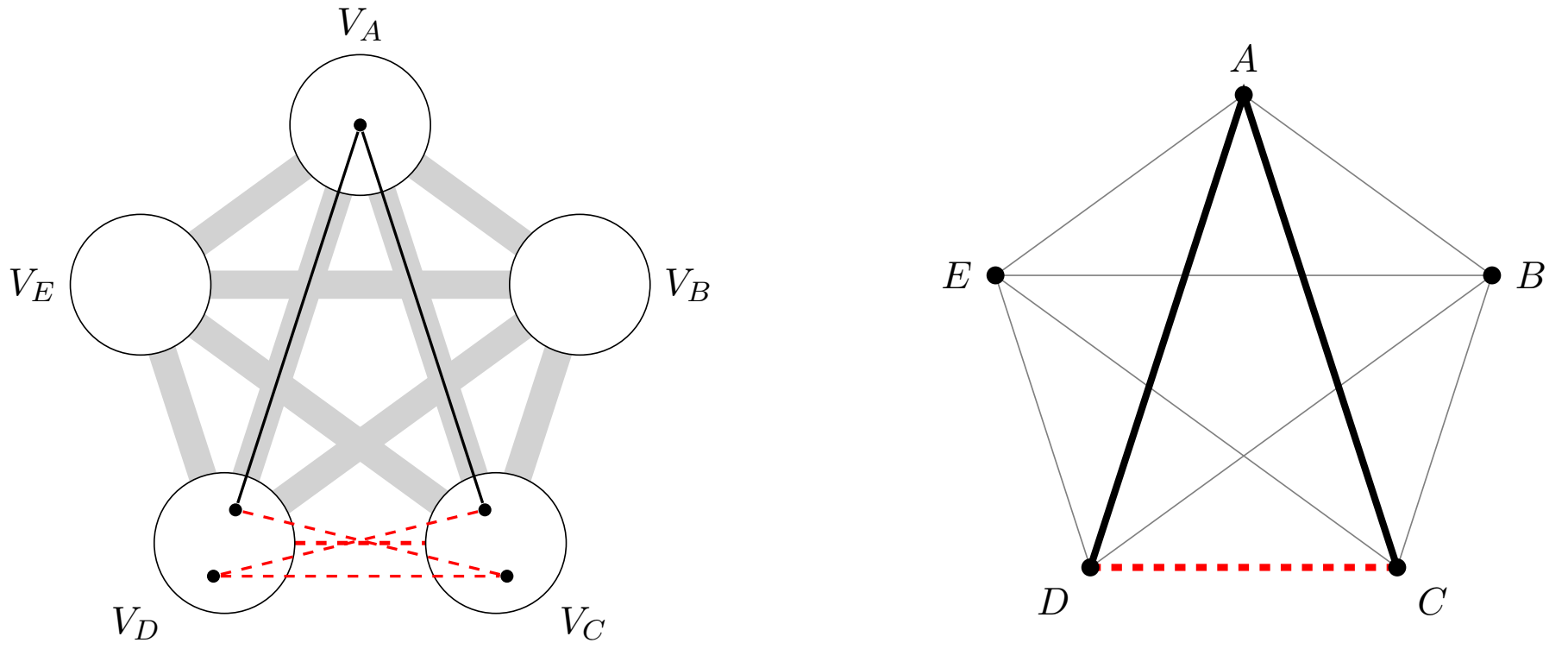
\begin{figure}[ht]
\centering
\scalebox{0.77}{%
\begin{tikzpicture}[
  part/.style={draw,circle,minimum size=1.4cm,inner sep=0pt,fill=white},
  band/.style={line width=8pt,gray!35,line cap=round},
  banddel/.style={line width=1.2pt,red,dashed},
  v/.style={fill,circle,inner sep=1.3pt},
  pt/.style={fill,circle,inner sep=1.8pt},
  present/.style={thick},
  missing/.style={thick,red,dashed}]

\begin{scope}
\coordinate (A) at (90:2.3);
\coordinate (B) at (18:2.3);
\coordinate (C) at (-54:2.3);
\coordinate (D) at (-126:2.3);
\coordinate (E) at (162:2.3);
\foreach \p/\q in {A/B,A/C,A/D,A/E,B/C,B/D,B/E,C/E,D/E}{\draw[band] (\p)--(\q);}
\draw[banddel] (C)--(D);
\node[part,label=above:$V_A$]        at (A) {};
\node[part,label=right:$V_B$]        at (B) {};
\node[part,label=below right:$V_C$]  at (C) {};
\node[part,label=below left:$V_D$]   at (D) {};
\node[part,label=left:$V_E$]         at (E) {};

\node[v] (x5) at (A) {};
\node[v] (x1) at ($(C)+(-0.11, 0.33)$) {};
\node[v] (x3) at ($(C)+( 0.11,-0.33)$) {};
\node[v] (x4) at ($(D)+( 0.11, 0.33)$) {};
\node[v] (x2) at ($(D)+(-0.11,-0.33)$) {};
\draw[present] (x4)--(x5) (x5)--(x1);
\draw[missing] (x1)--(x2) (x2)--(x3) (x3)--(x4);

\node at (0,-3.3) {\small (a) A copy of $C_5$ before the deletion.};
\end{scope}

\begin{scope}[shift={(8.8,0)}]
\coordinate (A) at (90:2.6);
\coordinate (B) at (18:2.6);
\coordinate (C) at (-54:2.6);
\coordinate (D) at (-126:2.6);
\coordinate (E) at (162:2.6);
\foreach \p/\q in {A/B,A/E,B/C,B/D,B/E,C/E,D/E}{\draw[gray] (\p)--(\q);}
\draw[line width=2pt] (D)--(A)--(C);
\draw[line width=2pt,red,dashed] (C)--(D);
\node[pt,label=above:$A$]       at (A) {};
\node[pt,label=right:$B$]       at (B) {};
\node[pt,label=below right:$C$] at (C) {};
\node[pt,label=below left:$D$]  at (D) {};
\node[pt,label=left:$E$]        at (E) {};
\node at (0,-3.3) {\small (b) The graph $\mathcal{T}$ and the projected copy.};
\end{scope}

\end{tikzpicture}
}
\caption{The construction for $r=2$, $d=1$ and $\mathcal{F}=C_5$ with
five parts. The five edges of the copy of $C_5$ correspond to only the
three edges of the triangle $ACD$, a smallest homomorphic image of
$C_5$, so deleting $CD$ from $K_5$ destroys the copy.}
\label{figure: construction}
\end{figure}

\paragraph{Proof sketch.} 
Note that a blow-up of an $rd$-graph $\mathcal{T}$ on $m$ vertices with
sufficiently large parts contains a copy of $\mathcal{F}^d$ if and only
if $\mathcal{F}^d\rightarrow\mathcal{T}$.
Moreover, if the parts have sizes as equal as possible, the density of
the blow-up is close to that of $\mathcal{T}$, since a random $rd$-set
meets some part twice with probability only about $\binom{rd}{2}/m$.
It therefore suffices to find a dense $rd$-graph $\mathcal{T}$ on $m$
vertices, where $m$ is independent of $n$, such that there is no
homomorphism from $\mathcal{F}^d$ to $\mathcal{T}$.
To this end, we consider $K^{(rd)}_m$ and delete each edge
independently with probability $p$, chosen later.
A fixed homomorphism $\mathcal{F}^d\rightarrow K^{(rd)}_m$ survives
the deletion only if no edge of its image is deleted.
By Lemma~\ref{lemma: edge_number_product}, its image has at least
$\ehom(\mathcal{F}^d)=\ehom(\mathcal{F})^d$ edges, so it survives with probability at most $(1-p)^{\ehom(\mathcal{F})^d}$.
Note that the exponent in this bound is the homomorphic edge number of
$\mathcal{F}^d$ and not its number of edges, since several edges of a
copy of $\mathcal{F}^d$ in the blow-up may correspond to the same edge
of $\mathcal{T}$, see Figure~\ref{figure: construction} for an example.

There are only $m^{|V(\mathcal{F})|d}$ maps from $V(\mathcal{F}^d)$ to
$[m]$, a number which is crucially independent of $n$.
The union bound therefore shows that, with positive probability,
no homomorphism survives provided that
$p\,\ehom(\mathcal{F})^d>|V(\mathcal{F})|\,d\log m$.
We therefore choose $p$ slightly larger than
$(|V(\mathcal{F})|\,d\log m)/\ehom(\mathcal{F})^d$.
The two losses, $p$ from the deletion and $\binom{rd}{2}/m$ from the
blow-up, are balanced by $m=\ehom(\mathcal{F})^d$, and both are then
of order $d^2\ehom(\mathcal{F})^{-d}$ up to constants depending only
on $\mathcal{F}$.
For powers, our construction is tight up to a polynomial factor in
$d$, since the homomorphic edge number is multiplicative by
Lemma~\ref{lemma: edge_number_product}.
Indeed, a smallest homomorphic image of $\mathcal{F}^d$ has
$\ehom(\mathcal{F})^d$ edges, and averaging over such an image gives
the matching lower bound.
We now give the formal proof.

\begin{proof}[Proof of Theorem~\ref{theorem: products}]
Let $\mathcal{F}$ be a fixed $r$-graph and let $\ehom=\ehom(\mathcal{F})$.
If $\mathcal{F}$ is $r$-partite, then $\mathcal{F}^d$ is $rd$-partite,
and thus $\pi(\mathcal{F}^d)=0$ by Lemma~\ref{lemma: standard}.
Suppose now that $\mathcal{F}$ is not $r$-partite, in which case
$\ehom\geq2$.

We begin with the lower bound $1-\pi(\mathcal{F}^d)\geq\ehom^{-d}$.
Let $\varphi_0$ be a homomorphism from $\mathcal{F}$ to an $r$-graph
$\mathcal{H}$ with $|E(\mathcal{H})|=\ehom$.
The map that applies $\varphi_0$ in each of the $d$ factors is a
homomorphism $\mathcal{F}^d\rightarrow\mathcal{H}^d$, so
$\pi(\mathcal{F}^d)\leq\pi(\mathcal{H}^d)$ by
Lemma~\ref{lemma: standard}.
Note that $|E(\mathcal{H}^d)|=|E(\mathcal{H})|^d=\ehom^d$.
Let $\mathcal{G}$ be an $\mathcal{H}^d$-free $rd$-graph on
$n\geq|V(\mathcal{H}^d)|$ vertices, and pick an injective map
$\psi:V(\mathcal{H}^d)\rightarrow V(\mathcal{G})$ uniformly at random.
For every $B\in E(\mathcal{H}^d)$, the set $\psi(B)$ is a uniformly
random $rd$-subset of $V(\mathcal{G})$, so $\psi(B)\notin E(\mathcal{G})$
with probability $1-|E(\mathcal{G})|/\binom{n}{rd}$.
Since $\mathcal{G}$ is $\mathcal{H}^d$-free, we have that for every choice of $\psi$ at least one $B\in E(\mathcal{H}^d)$ satisfies
$\psi(B)\notin E(\mathcal{G})$, so the expected number of such edges is
at least one, that is,
\[
    \ehom^d\left(1-\frac{|E(\mathcal{G})|}{\binom{n}{rd}}\right)\geq1.
\]
Hence $\ex(n,\mathcal{H}^d)\leq(1-\ehom^{-d})\binom{n}{rd}$ for every
$n\geq|V(\mathcal{H}^d)|$, and letting $n\rightarrow\infty$ gives
$1-\pi(\mathcal{F}^d)\geq1-\pi(\mathcal{H}^d)\geq\ehom^{-d}$, as claimed.

We next show that $1-\pi(\mathcal{F}^d)\leq Cd^2\,\ehom^{-d}$, following
the construction sketched above.
Fix $d\geq1$ and let $m=\ehom^d$ and 
\[
    p=\frac{|V(\mathcal{F})|(\log\ehom)\,d^2+d\log\ehom}{\ehom^d}.
\]
We shall prove that $1-\pi(\mathcal{F}^d)\leq p+\binom{rd}{2}/m+\ehom^{-d}$,
which implies the claimed bound.
If $p>1$ or $m<rd$, this holds trivially, as the right-hand side is
then larger than one, so we may assume that $p\leq1$ and $m\geq rd$.
Let $n\geq m$ and partition $[n]$ into parts $V_1,\ldots,V_m$, each of
size $\lfloor n/m\rfloor$ or $\lceil n/m\rceil$.
Let $\mathcal{K}$ be the complete $m$-partite $rd$-graph on $[n]$, that
is, the $rd$-graph whose edges are the $rd$-sets meeting $rd$ distinct
parts.
For an edge $\sigma$ of the complete $rd$-graph $K^{(rd)}_m$ on
$[m]$, let $\mathcal{K}_\sigma$ be the set of edges of $\mathcal{K}$
having exactly one vertex in $V_i$ for every $i\in\sigma$.
Then $E(\mathcal{K})$ is the disjoint union of the sets $\mathcal{K}_\sigma$
over all $\binom{m}{rd}$ edges $\sigma$ of $K^{(rd)}_m$, and we have
\begin{equation}\label{eq: complete multipartite density}
    |E(\mathcal{K})|
    \geq\binom{m}{rd}\left\lfloor\frac{n}{m}\right\rfloor^{rd}
    \geq(1-o(1))\left(1-\frac{\binom{rd}{2}}{m}\right)\binom{n}{rd},
\end{equation}
where the $o(1)$-term tends to zero as $n\rightarrow\infty$ and we used that $\left(\binom{m}{rd}(rd)!\right)/m^{rd}
=\prod_{k=0}^{rd-1}\left(1-k/m\right)
\geq1-\binom{rd}{2}/m$.
Now delete each edge of $K^{(rd)}_m$ independently with probability
$p$, and let $\mathcal{G}$ be the $rd$-graph on $[n]$ whose edge set is
the union of the sets $\mathcal{K}_\sigma$ over all edges $\sigma$ of
$K^{(rd)}_m$ that were not deleted.
In other words, if $\mathcal{T}$ denotes the $rd$-graph on $[m]$ given
by the edges that were not deleted, then $\mathcal{G}$ is the blow-up
of $\mathcal{T}$ with parts $V_1,\ldots,V_m$.
Let $\mathcal{A}$ be the event that every homomorphism
$\varphi:\mathcal{F}^d\rightarrow K^{(rd)}_m$ has at least one edge
$B\in E(\mathcal{F}^d)$ for which the edge $\varphi(B)$ of
$K^{(rd)}_m$ was deleted.

We claim that, if $\mathcal{A}$ occurs, then $\mathcal{G}$ is
$\mathcal{F}^d$-free.
Indeed, consider a copy of $\mathcal{F}^d$ in $\mathcal{G}$ and map
every vertex of $\mathcal{F}^d$ to the index of the part containing its
image in the copy.
Every edge of the copy is an edge of $\mathcal{K}$ and thus meets $rd$
distinct parts, so this map is a homomorphism
$\varphi:\mathcal{F}^d\rightarrow K^{(rd)}_m$, and the edge of the copy
corresponding to $B\in E(\mathcal{F}^d)$ is contained in
$\mathcal{K}_{\varphi(B)}$.
Since it is an edge of $\mathcal{G}$, the edge $\varphi(B)$ of
$K^{(rd)}_m$ was not deleted.
As this holds for every $B\in E(\mathcal{F}^d)$, the event $\mathcal{A}$ does not occur, a contradiction.

We next show that $\Pr(\compl{\mathcal{A}})\leq\ehom^{-d}$. 
To this end, fix a homomorphism $\varphi:\mathcal{F}^d\rightarrow K^{(rd)}_m$.
TThe $rd$-graph on $\varphi(V(\mathcal{F}^d))$ with edge set
$\{\varphi(B):B\in E(\mathcal{F}^d)\}$ is a homomorphic image of
$\mathcal{F}^d$ and thus has at least $\ehom^d$ edges by
Lemma~\ref{lemma: edge_number_product}.
The edges of $K^{(rd)}_m$ are deleted independently and thus the
probability that none of the edges of this image is deleted is at most
$(1-p)^{\ehom^d}\leq\exp(-p\,\ehom^d)$.
Since there are at most
$m^{|V(\mathcal{F})|d}=\exp\bigl(|V(\mathcal{F})|(\log\ehom)d^2\bigr)$
maps from $V(\mathcal{F}^d)$ to $[m]$, the union bound and the choice of
$p$ give
$\Pr(\compl{\mathcal{A}})\leq m^{|V(\mathcal{F})|d} \cdot \exp(-p\,\ehom^d)=\ehom^{-d}$.

Each edge of $K^{(rd)}_m$ is kept with probability $1-p$, so
$\mathbb{E}(|E(\mathcal{G})|)=(1-p)|E(\mathcal{K})|$.
Let $\mathbf{1}_{\mathcal A}$ denote the indicator of the event
$\mathcal A$, i.e., $|E(\mathcal G)|\,\mathbf{1}_{\mathcal A}$ equals
$|E(\mathcal G)|$ if $\mathcal A$ occurs and zero otherwise.
Since $|E(\mathcal{G})|\leq\binom{n}{rd}$, $\mathbb{E}\bigl(|E(\mathcal{G})|\,\mathbf{1}_{\mathcal{A}}\bigr) = \mathbb{E}(|E(\mathcal{G})|) - \mathbb{E}\bigl(|E(\mathcal{G})|\,\mathbf{1}_{\compl{\mathcal A}}\bigr)$, and
$\Pr(\compl{\mathcal{A}})\leq\ehom^{-d}$, we have
\[
    \mathbb{E}\bigl(|E(\mathcal{G})|\,\mathbf{1}_{\mathcal{A}}\bigr) \geq \mathbb{E}(|E(\mathcal{G})|) - \ehom^{-d}\binom{n}{rd}.
\]
Together with~\eqref{eq: complete multipartite density}, this gives
\[
    \mathbb{E}\bigl(|E(\mathcal{G})|\,\mathbf{1}_{\mathcal{A}}\bigr) \geq
    \left(1-p-\frac{\binom{rd}{2}}{m}-\ehom^{-d}-o(1)\right)\binom{n}{rd}.
\]
Since $\Pr(\mathcal{A})\geq1-\ehom^{-d}>0$ and
$|E(\mathcal G)|\,\mathbf{1}_{\mathcal A}$ is zero whenever
$\mathcal A$ does not occur, there is a choice of deleted edges of
$K^{(rd)}_m$ for which $\mathcal A$ occurs and $|E(\mathcal G)|$ is
at least the right-hand side of this inequality.
For this choice, $\mathcal G$ is $\mathcal F^d$-free.
Letting $n\rightarrow\infty$ gives
\[
    1-\pi(\mathcal{F}^d)
    \leq
    p+\frac{\binom{rd}{2}}{m}+\ehom^{-d},
\]
as claimed.
\end{proof}

\paragraph{Remark.}
The proof above gives, for all integers $d\geq1$ and $m\geq rd$, the
upper bound
\begin{equation}\label{eq: general bound}
    1-\pi(\mathcal{F}^d)
    \leq\frac{|V(\mathcal{F})|\,d\log m+d\log\ehom}{\ehom^d}
    +\frac{\binom{rd}{2}}{m}+\ehom^{-d},
\end{equation}
where, as in the proof, $\ehom=\ehom(\mathcal{F})$.
The choice $m=\ehom^d$ was only used to balance the first two terms.

\subsection{Proof of \texorpdfstring{Proposition~\ref{theorem: two factors}}{Proposition 8}}\label{section: two factors}

The idea of the proof is to first replace $\mathcal{F}$ by a smallest
homomorphic image $\mathcal{H}$, which does not decrease the Tur\'an
density of the product, and then to map the vertices of $\mathcal{H}$
injectively and at random into an $(\mathcal{H}*\mathcal{G})$-free
$(r+s)$-graph of maximum size.
The common link of the images of the edges of $\mathcal{H}$, that is,
the $s$-graph of all $s$-sets whose union with each of these images is
an edge, is then $\mathcal{G}$-free, and a straightforward averaging
argument yields the claimed bound.

\begin{proof}[Proof of Proposition~\ref{theorem: two factors}]
Since $\mathcal{F}*\mathcal{G}=\mathcal{G}*\mathcal{F}$, it suffices
to show that
$1-\pi(\mathcal{F}*\mathcal{G})\geq(1-\pi(\mathcal{G}))/\ehom(\mathcal{F})$.
Let $\varphi$ be a homomorphism from $\mathcal{F}$ to an $r$-graph
$\mathcal{H}$ with $|E(\mathcal{H})|=\ehom(\mathcal{F})$, where
$V(\mathcal{H})$ and $V(\mathcal{G})$ are disjoint.
The map $\varphi'$ on $V(\mathcal{F})\cup V(\mathcal{G})$ with
$\varphi'(x)=\varphi(x)$ for $x\in V(\mathcal{F})$ and $\varphi'(y)=y$
for $y\in V(\mathcal{G})$ is a homomorphism from
$\mathcal{F}*\mathcal{G}$ to $\mathcal{H}*\mathcal{G}$, as it maps
every edge $A\cup B$ of $\mathcal{F}*\mathcal{G}$ to the edge
$\varphi(A)\cup B$ of $\mathcal{H}*\mathcal{G}$.
Hence $\pi(\mathcal{F}*\mathcal{G})\leq\pi(\mathcal{H}*\mathcal{G})$ by
Lemma~\ref{lemma: standard}, and it suffices to show that
$1-\pi(\mathcal{H}*\mathcal{G})\geq(1-\pi(\mathcal{G}))/\ehom(\mathcal{F})$.

Let $\mathcal{K}$ be an $(\mathcal{H}*\mathcal{G})$-free $(r+s)$-graph
on $n\geq|V(\mathcal{H})|+s$ vertices with
$|E(\mathcal{K})|=\ex(n,\mathcal{H}*\mathcal{G})$,
and let $\varepsilon=\varepsilon(n)=|E(\mathcal{K})|/\binom{n}{r+s}$
denote its edge density, so
$\pi(\mathcal{H}*\mathcal{G})=\lim_{n\rightarrow\infty}\varepsilon$.
Choose now an injective map
$\psi:V(\mathcal{H})\rightarrow V(\mathcal{K})$ uniformly at random and
let $X=V(\mathcal{K})\setminus\psi(V(\mathcal{H}))$.
Let $\mathcal{L}$ be the $s$-graph on $X$ whose edges are the sets
$B\in\binom{X}{s}$ such that $\psi(A)\cup B$ is an edge of $\mathcal{K}$ for every
$A\in E(\mathcal{H})$.
In other words, $\mathcal{L}$ is the common link in $\mathcal{K}$ of the sets $\psi(A)$, where $A\in E(\mathcal{H})$.
Since any copy of $\mathcal{G}$ in $\mathcal{L}$, together with $\psi(\mathcal{H})$, forms a copy of $\mathcal{H}*\mathcal{G}$ in $\mathcal{K}$, the $s$-graph $\mathcal{L}$ is $\mathcal{G}$-free, and thus $|E(\mathcal{L})|\leq\ex(|X|,\mathcal{G})$. 

After choosing $\psi$, we pick $B$ uniformly at random from
$\binom{X}{s}$, so for every fixed $A\in E(\mathcal{H})$,
the set $\psi(A)\cup B$ is chosen uniformly at random from
$\binom{V(\mathcal{K})}{r+s}$.
Indeed, every possible pair $(\psi,B)$ is equally likely.
For a fixed $U\in\binom{V(\mathcal{K})}{r+s}$, the pairs satisfying
$\psi(A)\cup B=U$ are obtained by choosing an $r$-set
$R\in\binom{U}{r}$, a bijection from $A$ to $R$, and an injection
from $V(\mathcal{H})\setminus A$ to
$V(\mathcal{K})\setminus U$, after which $B=U\setminus R$.
Note that the number of such pairs does not depend on $U$.
Thus $\psi(A)\cup B$ is not an edge of $\mathcal{K}$ with probability
$1-\varepsilon$.
By the definition of $\mathcal{L}$, we have that $B\notin E(\mathcal{L})$
exactly when $\psi(A)\cup B\notin E(\mathcal{K})$ for some
$A\in E(\mathcal{H})$.
Hence the union bound gives
\[
    \mathbb E(|E(\mathcal{L})|)
    =\binom{|X|}{s} \Pr(B\in E(\mathcal{L})) \geq
    \binom{|X|}{s} \bigl(1-(1-\varepsilon) \ehom(\mathcal{F})\bigr).
\]
Since $|E(\mathcal{L})|\leq\ex(|X|,\mathcal{G})$ for every choice of
$\psi$, it follows that
\[
    \ex(|X|,\mathcal{G})
    \geq \binom{|X|}{s} \bigl(1-(1-\varepsilon) \ehom(\mathcal{F})\bigr).
\]
Dividing by $\binom{|X|}{s}$ and letting $n\rightarrow\infty$ gives
$\pi(\mathcal{G})\geq1-\bigl(1-\pi(\mathcal{H}*\mathcal{G})\bigr) \ehom(\mathcal{F})$,
and rearranging yields the claimed bound.
\end{proof}

\section{Concluding remarks}\label{section: remarks}

In this note, we determined the order of
$1-\lim_{r\rightarrow\infty}\pi(\mathcal{D}_r(s,t))$ for every fixed
$s\geq3$ as $t \rightarrow \infty$, and we showed that the homomorphic edge number is a useful parameter for Tur\'an problems on hypergraph products.
In terms of this parameter, we determined the exponential rate at
which $1-\pi(\mathcal{F}^d)$ tends to zero for every $r$-graph
$\mathcal{F}$ that is not $r$-partite and gave a lower bound on
$1-\pi(\mathcal{F}*\mathcal{G})$ for arbitrary $\mathcal{F}$ and
$\mathcal{G}$.
We now close this note with some remarks on the limitations of our results and suggest further questions.

\paragraph{Daisies.}
For every fixed $s\geq3$, we have determined the order of
$1-\lim_{r\rightarrow\infty}\pi(\mathcal{D}_r(s,t))$ as
$t\rightarrow\infty$.
The bound of de~Caen gives the lower bound
$(1+o(1))(s-1)!\,t^{-(s-1)}$, and taking $q$ in the proof of
Theorem~\ref{main} as the largest prime with $q\leq(t-s)/(s-1)$,
rather than a power of two, gives the upper bound
$(1+o(1))(s-1)^{s-1}t^{-(s-1)}$. 
Here, we use the prime number theorem, which gives
$q\sim t/(s-1)$ for fixed $s$ as $t\rightarrow\infty$.
It would be interesting to determine whether a leading
constant exists and, if so, to determine its value.

\paragraph{Hypergraph powers.}
For non-$r$-partite $\mathcal{F}$,
Theorem~\ref{theorem: products} determines the exponential rate of
$1-\pi(\mathcal{F}^d)$ as $d\rightarrow\infty$.
For $d=1$, the lower bound $1-\pi(\mathcal{F})\geq\ehom(\mathcal{F})^{-1}$
is the standard averaging bound applied to a smallest homomorphic image
of $\mathcal{F}$, while the upper bound contains an unspecified
constant depending on $\mathcal{F}$.
We also note that the lower bound can be sharpened by directly applying
Proposition~\ref{theorem: two factors}.
Indeed, applying it to $\mathcal{F}^{d-1}*\mathcal{F}$ and using
Lemma~\ref{lemma: edge_number_product} gives $1-\pi(\mathcal{F}^d) \geq (1-\pi(\mathcal{F}))/\ehom(\mathcal{F})^{d-1}$ for every $d\geq2$, and the same inequality holds with equality
for $d=1$.
Moreover, the bound of Sidorenko~\cite{Sid}
(see also~\cite{Sid94}), applied to a smallest homomorphic image
of $\mathcal{F}$, gives $1-\pi(\mathcal{F})\geq1/(\ehom(\mathcal{F})-1)$, and thus
\[
    1-\pi(\mathcal{F}^d)
    \geq
    \frac{1}{(\ehom(\mathcal{F})-1)\ehom(\mathcal{F})^{d-1}}.
\]
Since these refinements do not change the dependence on $d$,
we have opted to present the simpler lower bound in
Theorem~\ref{theorem: products} to emphasize that its upper and
lower bounds agree up to a factor of order $d^2$.
In any case, the bounds do not generally determine the exact value of
$\pi(\mathcal{F}^d)$ for a fixed $d$.
We leave it as an open question whether the factor $d^2$ can be
omitted, that is, whether
$1-\pi(\mathcal{F}^d)= \Theta_{\mathcal{F}}\bigl(\ehom(\mathcal{F})^{-d}\bigr)$ as $d\rightarrow\infty$.

\paragraph{Products of two different factors.}
Proposition~\ref{theorem: two factors} only gives a lower bound on
$1-\pi(\mathcal{F}*\mathcal{G})$, and using the
construction of Section~\ref{section: powers} does not in general
determine $\pi(\mathcal{F}*\mathcal{G})$.
The following example shows, however, that
$\pi(\mathcal{F}*\mathcal{G})$ is not determined by the Tur\'an
densities and the uniformities of the factors alone.

For $k\geq2$, let $\mathcal{A}_k$ be the $2k$-graph with three edges
$E_1\cup E_2$, $E_1\cup E_3$, $E_2\cup E_3$ for disjoint $k$-sets
$E_1,E_2,E_3$, and let $\mathcal{B}_k$ be the disjoint union of
$\mathcal{A}_k$ and a $2k$-graph $\mathcal{C}_k$ consisting of $k+1$
distinct $2k$-subsets of a $(2k+1)$-set.
By a result of Frankl~\cite{Frankl}, $\pi(\mathcal{A}_k)=1/2$.
Moreover, a $\mathcal{C}_k$-free $2k$-graph on $n$ vertices has at most
$k$ edges inside every $(2k+1)$-set, and averaging over all
$(2k+1)$-sets gives $\pi(\mathcal{C}_k)\leq k/(2k+1)<1/2$.
Since the Tur\'an density of a vertex-disjoint union of two
hypergraphs is the maximum of their Tur\'an densities, we obtain
$\pi(\mathcal{B}_k)=1/2$.
On the other hand, $\ehom(\mathcal{A}_k)=3$ and
$\ehom(\mathcal{B}_k)\geq k+1$, since in both $\mathcal{A}_k$ and
$\mathcal{C}_k$ every two vertices are contained in a common edge.
Proposition~\ref{theorem: two factors} therefore gives
$1-\pi(\mathcal{A}_k^2)\geq(1-\pi(\mathcal{A}_k))/\ehom(\mathcal{A}_k)=1/6$,
while~\eqref{eq: general bound}, with $m=k^3$, $|V(\mathcal{B}_k)|=5k+1$ and $\ehom(\mathcal{B}_k)\geq k+1$, gives
$1-\pi(\mathcal{B}_k^2)=O(\log k/k)$, so $\pi(\mathcal{B}_k^2)\rightarrow1$
as $k\rightarrow\infty$.
Hence, if $\pi(\mathcal{F}*\mathcal{G})$ were determined by
$\pi(\mathcal{F})$, $\pi(\mathcal{G})$ and the uniformities of
$\mathcal{F}$ and $\mathcal{G}$, then $\pi(\mathcal{A}_k^2)$ and
$\pi(\mathcal{B}_k^2)$ would be equal, since $\mathcal{A}_k$ and
$\mathcal{B}_k$ are $2k$-graphs with Tur\'an density $1/2$.
However, $\pi(\mathcal{A}_k^2)\leq5/6$, while
$\pi(\mathcal{B}_k^2)>5/6$ for all sufficiently large $k$, and this
difference is explained by the homomorphic edge numbers of
$\mathcal{A}_k$ and $\mathcal{B}_k$.

\paragraph{Ramsey numbers of daisies and hypergraph products.}
For an $r$-graph $\mathcal{H}$ and an integer $\ell\geq2$,
let $R(\mathcal{H};\ell)$ be the least $n$ such that every
$\ell$-colouring of $E(K_n^{(r)})$ contains a monochromatic
copy of $\mathcal{H}$.
Pudl\'ak, R\"odl and Sales~\cite{PRS} studied Ramsey
numbers of daisies and obtained bounds for both fixed and
unbounded stem sizes.
Analogously to Problems~\ref{problem: two_product}
and~\ref{problem: fold}, they also asked how
$R(\mathcal{F}*\mathcal{G};\ell)$ compares to
$R(\mathcal{F};\ell)$ and $R(\mathcal{G};\ell)$, and,
more specifically, how $R((K_3^{(2)})^d;\ell)$ grows as
a function of $d$ and $\ell$,
see~\cite[Problems~7.4 and~7.5]{PRS}.
Perhaps our methods can also be used to obtain bounds for
Ramsey numbers of daisies and more general hypergraph products.

\paragraph{AI declaration.}
Theorem~\ref{main} grew out of the author's earlier joint work with Axenovich on Tur\'an densities of subgraphs of the hypercube, see~\cite[Proposition~2.5]{AP}, where ideas from the construction of Ellis, Ivan and Leader~\cite{EIL} were used independently of any application to daisies.
Based on this, the author suggested to GPT-6 Astra the specific
modification of the construction presented in the proof of
Theorem~\ref{main}. 
A subsequent brief iterative exchange, guided by
the author, ultimately led to a proof of the theorem.
For Theorem~\ref{theorem: products} and Proposition~\ref{theorem: two factors},
GPT-6 Astra was used to help work out technical details of the
proofs. 
In both cases, all arguments in which GPT-6 Astra assisted were
checked, shortened and reworked by the author.
GPT-6 Astra was also used to assist with the literature search on
Tur\'an problems for hypergraphs.
Finally, the note was written by the author, with GPT-6 Astra used to find typos and misprints, and occasionally to check grammar.

\end{document}